\documentclass[10pt,a4paper,oneside,leqno]{article}
\usepackage{amsfonts,amssymb}
\usepackage{amscd}
\usepackage{color}
\newtheorem{defn}{Definition}
\newtheorem{conven}{Convention}

\newtheorem{observen}{Observation}

\vspace{0.5cm}
 4
\font\ebf=cmbx8
\font\erm=cmr8

\usepackage{graphicx, graphics}
\usepackage{color}

\newcommand{\os}{\oplus\!\!\to}
\newcommand{\cs}{\copyright\!\!\to}

\begin{document}
\begin{center}
	\noindent { \textsc{ Graded posets  zeta matrix  formula}}  \\ 
	\vspace{0.3cm}
	\vspace{0.3cm}
	\noindent Andrzej Krzysztof Kwa\'sniewski \\
	\vspace{0.2cm}
	\noindent {\erm Member of the Institute of Combinatorics and its Applications  }\\
{\erm High School of Mathematics and Applied Informatics} \\
	{\erm  Kamienna 17, PL-15-021 Bia\l ystok, Poland }\noindent\\
	\noindent {\erm e-mail: kwandr@gmail.com}\\
	\vspace{0.4cm}
\end{center}

\noindent {\ebf Abstract:}
\vspace{0.1cm}
\noindent {\small The way to arrive at formula  of zeta matrix for any graded posets  with the finite set of minimal elements is delivered following [1].  This is being achieved via adjacency and zeta matrix description of bipartite digraphs chains  - the representatives of graded posets. The bipartite digraphs elements of such chains amalgamate to form corresponding cover relation graded poset digraphs with corresponding adjacency matrices being amalgamated  throughout natural join as special adequate database operation. The colligation  of  reachability and connectivity with the presented description  is made explicit. 
\noindent The special posets  encoded via  KoDAGs directed  acyclic graphs as cobeb posets Hasse diagrams  are  recognized as an example of  differential  posets subfamily.  
\noindent As on the    01.01.2009 one reminisce  261-th  anniversary  of death of  Johann Bernoulli the First this Sylvester Night article is to commemorate  this date.}

\vspace{0.3cm}

\noindent Key Words: graded digraphs, differential posets,  natural join, umbral calculus

\vspace{0.1cm}

\noindent AMS Classification Numbers: 06A06 ,05B20, 05C75 , 05A30

\vspace{0.1cm}

\noindent  affiliated to The Internet Gian-Carlo Rota Polish Seminar:

\noindent \emph{http://ii.uwb.edu.pl/akk/sem/sem\_rota.htm}

\section{Preliminaries: notation and terminology.} 
\textbf{1.1.} We shall try to keep track of   NIST  Dictionary of Algorithms and Data Structures terminology.
Abbreviation:  directed acyclic graph = DAG.
\vspace{0.1cm} 
\noindent Note:   The transitive closure of a directed acyclic graph or DAG is the reachability relation of the DAG and a strict partial order.

\noindent \textbf{1.2.}    The following convention scheme is adopted:
directed graph representatives of binary relations scheme  are:

\vspace{0.1cm}

\noindent 
\begin{center}
bipartite digraph representative $D_R =(A \times A,R) \leftrightarrow  R \subseteq  A\times A $ \\
$\equiv$  "just"   digraph representative
$ D(R)\equiv D_R=(A,E) ,  E \leftrightarrow R \subseteq A \times A ,$
bipartite digraph representative $D_R=(A \times B,R) \leftrightarrow  R \subseteq A \times B$
\end{center}

\vspace{0.1cm}

\noindent \textbf{1.3}    A directed path, is an oriented simple path with  all arcs of  the same direction  i.e.  all internal nodes have in- and out-degrees  equal one.

\vspace{0.1cm}

\noindent \textbf{Comment 1.}

\noindent "A  directed path is  a natural join of arcs that thus form a chain of vertices",  "A  chain of  coded data objects is  a natural join of their subsequent pairs".

\vspace{0.1cm}

\noindent Anticipated: coded data objects = relations (with varying arity allowed), binary relations, bipartite digraphs, adjacency matrices  (of graphs or digraphs):  [1], [2].

\vspace{0.1cm}

\noindent \textbf{Comment 2.}

\noindent Because of immense number of  applications of digraphs beyond mathematics - frequently successfully done  by non-mathematicians - it happens sometimes  that various names are given for the same notions and different names for the same objects. Let us recall and/or establish  some of them.
\vspace{0.1cm}

\noindent \textbf{Recall} 

\noindent Reachability is the ability reach  some other vertex from a given vertex in a directed graph. For a directed graph 
$D = ( \Phi, E), E\subseteq \Phi \times \Phi$  the reachability relation of $D$ is its transitive closure of $E$, i.e. the set of all ordered pairs $(s, t)$ of vertices in $\Phi$  for which there exist vertices $v_0 = s, v_1, . . . , v_e = t$ such that $(v_{i - 1}, v_i ) \in E$  for all $1 \leq i \leq e$.\\
We define here ( for directed graphs - more than nontrivial) the reachability = connectivity partial order relation $R$ over the nodes of the DAG as such that $xRy$  iff there exists  a directed path from $x$  to $y$ .

\vspace{0.2cm}

\noindent \textbf{1.4.}  (relations set sum) 

\noindent 
\begin{center}
reachability $\cup$ reflexibility   $=$  reflexive reachability, \\  reachability $=$ connectivity,
\end{center}

\vspace{0.1cm}

\noindent This means that:

\vspace{0.2cm}

\noindent \textbf{1.4.1} The reachability = connectivity relation is ${\bf R}^{\infty}={\displaystyle\bigcup_{k>0}}{R}^k
=  \:{\bf transitive\: closure}$ of ${\bf R}$, i.e.
$${\bf R}^{\infty}={ R}^1\cup{R}^2\cup\ldots\cup{ R}^n\cup\ldots\;\Leftrightarrow\;A(R^{\infty})=A({ R})^{\copyright 1}\vee
A({R})^{\copyright 2}\vee \ldots \vee A( R)^{\copyright n}\vee \ldots,$$
where $A({\bf R})$ is the Boolean  adjacency matrix of the
relation ${\bf R}$ simple digraph  and  $\copyright$ stays for Boolean product.

\vspace{0.2cm}
\noindent  The \textbf{reflexive}  reachability  relation $\zeta({\bf R})\equiv {\bf
R}^{*}$ is defined as \\

$${\bf R^*}=R^0\cup R^1\cup R^2\cup\ldots\cup R^n\cup\ldots  
\bigcup_{k\geq 0}{R}^k={\bf R}^{\infty}\cup {\bf I}_A =$$
 
\begin{center}
=  {\bf transitive} and {\bf reflexive} closure of $\bf{R}$  $\Leftrightarrow$
\end{center}
 
$$
\Leftrightarrow\;  A(R^{\infty})=A({ R})^{\copyright 0}\vee A({
R})^{\copyright 1}\vee A({ R})^{\copyright 2}\vee \ldots \vee A(R)^{\copyright n}\vee \ldots .
$$

\vspace{0.1cm}

\noindent\textbf{Comment 3.}   Colligate and identify $\zeta({\bf R})\equiv {\bf
R}^{*}$ with incidence algebra zeta function and with zeta matrix of the poset
associated to its Hasse digraph.

\vspace{0.2cm}
\noindent \textbf{1.4.2. } (Notation continued) \\

\noindent In what follows we shall use mathematical  terms:  reachability  and reflexive reachability  according to: 
put $R =  \prec\cdot$  which is Hasse diagram i.e. cover relation digraph - notation.  Then note the \textbf{schemes }below.

\vspace{0.2cm}
\noindent The partial order $\leq$ for locally finite poset $\Pi=\left\langle \Phi,\leq
 \right\rangle$ with respect to $\Pi$'s cover relation
 $\prec\cdot$ is:
 \begin{enumerate}
 \item $<$ {\bf the connectivity relation in Hasse digraph   i.e.
 digraph $D=\left\langle \Phi,\prec\cdot
 \right\rangle$.}
 \item $\leq$  {\bf is the reflexive reachability relation in Hasse digraph   i.e.   digraph $D=\left\langle \Phi,\prec\cdot
 \right\rangle$.  }
 \end{enumerate}

 \noindent {\bf Schemes:}
 $$ <=\prec\cdot ^{\infty}=\;connectivity\;of\;\prec\cdot$$
 $$ \leq=\prec\cdot ^{*}=\;reflexive \; reachability\;of\;\prec\cdot$$
 $$\prec\cdot ^{*}=\zeta(\prec\cdot).$$

\vspace{0.2cm}
\noindent \textbf{1.4.3.}     $\zeta(\prec\cdot)$  is also used to denote  zeta matrix of the graded poset $P(D)=\Pi=\left\langle \Phi,\leq
 \right\rangle$   associated  to $D=\left\langle \Phi,\prec\cdot
 \right\rangle$  [1,2]  which is equivalent to say that $P(D)= \left\langle \Phi,\leq
 \right\rangle$  =  transitive, reflexive  closure  of $D=\left\langle \Phi,\prec\cdot
 \right\rangle$.

\vspace{0.2cm}
\noindent  \textbf{1.5.}\\
\noindent In order to get\textbf{ complete} graded digraph connect any two vertices lying on consecutive levels with an arc keeping one direction -  say - upwards ( see  KoDAG  in [1,2]).

\section{Natural join of adjacency matrices }  
      
\vspace{0.1cm}
\noindent  \textbf{2.1.}  Locally finite poset $\Pi$  if fixed, is denominated by all its covering pairs and vice versa of course i.e. 

\begin{center}
 $  \Pi = \left\langle \Phi,\leq \right\rangle \      \Leftrightarrow \   \left\langle \Phi,\prec\cdot \right\rangle, $
\end{center}

\noindent as all properties of    order origin follow from those of  transitivity requirement - and vice versa of course.  Specifically recall-note the equivalence of descriptions:

\begin{center}
The complete graded poset  $\Leftrightarrow$     The complete graded digraph.
\end{center}

\noindent The arcs of any digraph $G = \left\langle \Phi,E \right\rangle $    with no multiple edges stay  automatically  for arcs of cover relation    $\prec\cdot$ in the corresponding  poset   $ \Pi = \left\langle \Phi,\leq \right\rangle$  for which the digraph $G = \left\langle \Phi,E \right\rangle $  becomes Hasse diagram  i.e.  we have:

\begin{center}
 $ \leq \:= \:\prec\cdot_* $ = reflexive reachability of $\prec\cdot$ \\
 
$ \prec\cdot_* \equiv {(I - \prec\cdot)}^{-1} \! \equiv \! \zeta(\prec\cdot)$  
\end{center}

\vspace{0.1cm}

\noindent Unless differently stated,  we  shall identify a digraph with  its adjacency matrix. In our context these are to be Hasse digraphs 
$D=\left\langle \Phi,\prec\cdot \right\rangle$ of  graded posets  $ \Pi = \left\langle \Phi,\leq \right\rangle$.

\vspace{0.2cm}

\noindent There are three standard widely used prevalent encodings, three ways of portraying partially ordered sets  $P(D) = \left\langle \Phi,\leq \right\rangle$  : Hasse diagrams $D=\left\langle \Phi,\prec\cdot \right\rangle$; zeta matrices   $\zeta({\bf \leq})$; and  cover matrices  $\zeta({\bf \prec\cdot }$).  These matrices are of course the adjacency matrices of the corresponding digraphs  $\left\langle \Phi,\leq \right\rangle$  and $\left\langle \Phi,\prec\cdot \right\rangle$.  In the  incidence algebra description of locally finite posets   $\zeta({\bf \leq})$ may be identified with the incidence function i.e. the characteristic function of the partial order $\leq$, (see [1] and references therein for the source cobweb posets examples of these objects). Here down we shall use adjacency and biadjacency nomenclature [1,2].

\vspace{0.2cm}

\noindent \textbf{Examples of} $\zeta({\bf \leq})$\\                     

\noindent Let  $F$ denotes arbitrary natural numbers valued sequence.  Let  $A_N$ be the Hasse matrix .i.e. adjacency matrix of cover relation $\prec\cdot$ digraph denominated by sequence $N$ [1]. Then the zeta matrix  $\zeta = (1-\mathbf{A}_N)^{-1\copyright}$ for the denominated by  $F = N$ cobweb poset is of the form [1]  
\vspace{1mm}
$$ \left[\begin{array}{ccccccccccccccccc}
1 & 1 & 1 & 1 & 1 & 1 & 1 & 1 & 1 & 1 & 1 & 1 & 1 & 1 & 1 & 1 & \cdots\\
0 & 1 & \textbf{\textcolor{blue}{0}} & 1 & 1 & 1 & 1 & 1 & 1 & 1 & 1 & 1 & 1 & 1 & 1 & 1 & \cdots\\
0 & 0 & 1 & 1 & 1 & 1 & 1 & 1 & 1 & 1 & 1 & 1 & 1 & 1 & 1 & 1 & \cdots\\
0 & 0 & 0 & 1 & \textbf{\textcolor{blue}{0}} & \textbf{\textcolor{blue}{0}} & 1 & 1 & 1 & 1 & 1 & 1 & 1 & 1 & 1 & 1 & \cdots\\
0 & 0 & 0 & 0 & 1 & \textbf{\textcolor{blue}{0}} & 1 & 1 & 1 & 1 & 1 & 1 & 1 & 1 & 1 & 1 & \cdots\\
0 & 0 & 0 & 0 & 0 & 1 & 1 & 1 & 1 & 1 & 1 & 1 & 1 & 1 & 1 & 1 & \cdots\\
0 & 0 & 0 & 0 & 0 & 0 & 1 & \textbf{\textcolor{blue}{0}} & \textbf{\textcolor{blue}{0}} & \textbf{\textcolor{blue}{0}} & 1 & 1 & 1 & 1 & 1 & 1 & \cdots\\
0 & 0 & 0 & 0 & 0 & 0 & 0 & 1 & \textbf{\textcolor{blue}{0}} & \textbf{\textcolor{blue}{0}} & 1 & 1 & 1 & 1 & 1 & 1 & \cdots\\
0 & 0 & 0 & 0 & 0 & 0 & 0 & 0 & 1 & \textbf{\textcolor{blue}{0}} & 1 & 1 & 1& 1 & 1 & 1 & \cdots\\
0 & 0 & 0 & 0 & 0 & 0 & 0 & 0 & 0 & 1 & 1 & 1 & 1 & 1 & 1 & 1 & \cdots\\
0 & 0 & 0 & 0 & 0 & 0 & 0 & 0 & 0 & 0 & 1 & \textbf{\textcolor{blue}{0}} & \textbf{\textcolor{blue}{0}} & \textbf{\textcolor{blue}{0}} & \textbf{\textcolor{blue}{0}} & 1 & \cdots\\
0 & 0 & 0 & 0 & 0 & 0 & 0 & 0 & 0 & 0 & 0 & 1 & \textbf{\textcolor{blue}{0}} & \textbf{\textcolor{blue}{0}} & \textbf{\textcolor{blue}{0}} & 1 & \cdots\\
0 & 0 & 0 & 0 & 0 & 0 & 0 & 0 & 0 & 0 & 0 & 0 & 1 & \textbf{\textcolor{blue}{0}} & \textbf{\textcolor{blue}{0}} & 1 & \cdots\\
0 & 0 & 0 & 0 & 0 & 0 & 0 & 0 & 0 & 0 & 0 & 0 & 0 & 1 & \textbf{\textcolor{blue}{0}} & 1 & \cdots\\
0 & 0 & 0 & 0 & 0 & 0 & 0 & 0 & 0 & 0 & 0 & 0 & 0 & 0 & 1 & 1 & \cdots\\
0 & 0 & 0 & 0 & 0 & 0 & 0 & 0 & 0 & 0 & 0 & 0 & 0 & 0 & 0 & 1 & \cdots\\
. & . & . & . & . & . & . & . & . & . & . & . & . & . & . & . & . \cdots\\
 \end{array}\right]$$
\vspace{1mm}   
\noindent \textbf{Figure $\zeta_N$.  The incidence matrix
$\zeta$ for the  natural numbers  i.e. $N$ - cobweb poset}

\vspace{0.1cm}
\noindent Note that the  matrix $\zeta$ representing uniquely its corresponding  cobweb poset does  exhibits  a staircase structure of zeros above the diagonal (see above, see below) which is characteristic to and characteristics of Hasse diagrams of \textbf{all}\textbf{ cobweb posets} while for graded posets it is characteristic too - this time all together with additional zeros right to the staircase - zeros generated by the bi-adjacency matrices $\left\langle B_k\right\rangle_{k\geq 0}$ chain as described in Observation 4. . Another wards:  it is the natural join $\os$ chain 
$$ 
\os_{k=0}^n B_k,\;\;n \in N \cup \{\infty\}
$$
which is equivalent to $\zeta$ function \textbf{characteristics} of any fixed $F$-denominated graded poset $(P,\leq)$. 

\vspace{1mm}
$$ \left[\begin{array}{ccccccccccccccccc}
1 & 1 & 1 & 1 & 1 & 1 & 1 & 1 & 1 & 1 & 1 & 1 & 1 & 1 & 1 & 1 & \cdots\\
0 & 1 & 1 & 1 & 1 & 1 & 1 & 1 & 1 & 1 & 1 & 1 & 1 & 1 & 1 & 1 & \cdots\\
0 & 0 & 1 & 1 & 1 & 1 & 1 & 1 & 1 & 1 & 1 & 1 & 1 & 1 & 1 & 1 & \cdots\\
0 & 0 & 0 & 1 & \textbf{\textcolor{red}{0}} & 1 & 1 & 1 & 1 & 1 & 1 & 1 & 1 & 1 & 1 & 1 & \cdots\\
0 & 0 & 0 & 0 & 1 & 1 & 1 & 1 & 1 & 1 & 1 & 1 & 1 & 1 & 1 & 1 & \cdots\\
0 & 0 & 0 & 0 & 0 & 1 & \textbf{\textcolor{red}{0}} & \textbf{\textcolor{red}{0}} & 1 & 1 & 1 & 1 & 1 & 1 & 1 & 1 & \cdots\\
0 & 0 & 0 & 0 & 0 & 0 & 1 & \textbf{\textcolor{red}{0}} & 1 & 1 & 1 & 1 & 1 & 1 & 1 & 1 & \cdots\\
0 & 0 & 0 & 0 & 0 & 0 & 0 & 1 & 1 & 1 & 1 & 1 & 1 & 1 & 1 & 1 & \cdots\\
0 & 0 & 0 & 0 & 0 & 0 & 0 & 0 & 1 & \textbf{\textcolor{red}{0}} & \textbf{\textcolor{red}{0}} & \textbf{\textcolor{red}{0}} & \textbf{\textcolor{red}{0}} & 1 & 1 & 1 & \cdots\\
0 & 0 & 0 & 0 & 0 & 0 & 0 & 0 & 0 & 1 & \textbf{\textcolor{red}{0}} & \textbf{\textcolor{red}{0}} & \textbf{\textcolor{red}{0}} & 1 & 1 & 1 & \cdots\\
0 & 0 & 0 & 0 & 0 & 0 & 0 & 0 & 0 & 0 & 1 & \textbf{\textcolor{red}{0}} & \textbf{\textcolor{red}{0}} & 1 & 1 & 1 & \cdots\\
0 & 0 & 0 & 0 & 0 & 0 & 0 & 0 & 0 & 0 & 0 & 1 & \textbf{0 }& 1 & 1 & 1 & \cdots\\
0 & 0 & 0 & 0 & 0 & 0 & 0 & 0 & 0 & 0 & 0 & 0 & 1 & 1 & 1 & 1 & \cdots\\
0 & 0 & 0 & 0 & 0 & 0 & 0 & 0 & 0 & 0 & 0 & 0 & 0 & 1 & \textbf{\textcolor{red}{0}} & \textbf{\textcolor{red}{0}} & \cdots\\
0 & 0 & 0 & 0 & 0 & 0 & 0 & 0 & 0 & 0 & 0 & 0 & 0 & 0 & 1 & \textbf{\textcolor{red}{0}} & \cdots\\
0 & 0 & 0 & 0 & 0 & 0 & 0 & 0 & 0 & 0 & 0 & 0 & 0 & 0 & 0 & 1 & \cdots\\
. & . & . & . & . & . & . & . & . & . & . & . & . & . & . & . & . \cdots\\
 \end{array}\right]$$

\vspace{1mm} \noindent \textbf{Figure $\zeta_F$.  The matrix
$\zeta$ for the Fibonacci cobweb poset associated to \textbf{$F$-KoDAG} Hasse digraph }

\vspace{0.2cm}
\noindent \textbf{Comment 4.} 
 The given  $F$-denominated staircase zeros structure above the diagonal of zeta matrix $zeta$ is the \textbf{unique characteristics} of  its corresponding  \textbf{$F$-KoDAG} Hasse digraphs.

\vspace{0.4cm}

\noindent \textbf{2.2.} \textbf{The natural join condition}.

\vspace{0.1cm}

\noindent The natural join operation is a binary operation  like $\Theta$  \textbf{operator in computer science} denoted here by $\os$ symbol  deliberately referring - in a quite  reminiscent manner - to direct sum  $\oplus$  of adjacency Boolean  matrices and - as matter of fact and  in effect - to direct the sum  $\oplus$  of corresponding biadjacency [reduced]  matrices of digraphs under natural join.

\noindent  $\os$ is a natural operator for sequences construction.  $\os$ operates on multi-ary relations according to the scheme:   $(n+k)_{ary}  \os (k+m)_{ary}$  =   $(n+ k +m)_{ary}$

\noindent For example:  $(1+1)_{ary} \os(1+1)_{ary} = (1+ 1 +1)_{ary}$ , binary $\os$ binary = ternary.

\vspace{0.2cm} 
\noindent Accordingly an action of  $\os$  on these multi-ary relations' digraphs adjacency matrices is to be designed soon in what follows.

\vspace{0.2cm}
\noindent \textbf{Domain-Codomain $F$-sequence condition} $\mathrm{dom}(R_{k+1}) = \mathrm{ran} (R_k)$,  $k=0,1,2,...$ [1].  

\vspace{0.2cm}
\noindent Consider any  natural number valued sequence $F = \{F_n\}_{n\geq 0}$. Consider then any  chain of binary relations  defined on pairwise disjoint finite sets with cardinalities appointed by  $F$ -sequence elements values. For that to start we specify at first   a relations'  domain-co-domain $F$ - sequence.   

\vspace{0.2cm}
\noindent \textbf{Domain-Codomain  $F$-sequence $(|\Phi_n| = F_n )$}

$$
	\Phi_0,\Phi_1,...\Phi_i,...\ \ \Phi_k\cap\Phi_n = \emptyset \ \ for \ \ k \neq n, |\Phi_n|=F_n; \ \ i,k,n=0,1,2,...
$$

\noindent Let $\Phi=\bigcup_{k=0}^n\Phi_k$  be the corresponding ordered partition  [ anticipating -  $\Phi$ is the vertex set of  $D = (\Phi,\prec\!\!\cdot$ )   and  its transitive, reflexive closure $(\Phi, \leq)$] .     Impose $\mathrm{dom} (R_{k+1}) = \mathrm{ran} (R_k)$ condition , $k\in N \cup \{\infty\}$. What we get is binary relations chain.

\begin{defn} (Relation`s chain)
\noindent Let   $\Phi=\bigcup_{k=0}^n\Phi_k$ ,   $\Phi_k  \cap \Phi_n = \emptyset$  for  $k  \neq  n$ be  the ordered partition of the set  $\Phi$ .
\vspace{0.1cm}
\noindent Let a sequence of binary relations be given such that 
$$
	R_0,R_1,...,R_i,...,R_{i+n},...,\ \ R_k\subseteq\Phi_k\times\Phi_{k+1},\ \ \mathrm{dom}(R_{k+1}) = \mathrm{ran}(R_k).
$$

\noindent Then the sequence  $\langle R_k\rangle_{k\geq 0}$    is called natural join  (binary) \textbf{relation's chain}.
\end{defn}

\vspace{0.2cm}
\noindent Extension to varying arity relations' natural join chains is straightforward.

\vspace{0.2cm}
\noindent As necessarily $\mathrm{dom}(R_{k+1}) = \mathrm{ran}(R_k)$  for relations' natural join chain any given  binary relation's chain is not just a sequence  therefore we use "link to link " notation for $k, i , n = 1,2,3,...$ ready for  relational data basis applications: 
$$
 R_0 \os R_1 \os ... \os R_i \os ... \os R_{i+n},... is\ an\  F-chain\ of\  binary\ relations
$$

\noindent where  $\os$  denotes natural join of relations as well as both  natural join of their bipartite digraphs and the natural join of their representative adjacency matrices (see [1,2]).

\vspace{0.2cm}
\noindent Relation's $F$-chain  naturally  represented by [identified with] the chain of theirs  \textbf{bipartite digraphs}

$$
	{ R_0 \os R_1 \os ... \os R_i \os ... \os R_{i+n},... \Leftrightarrow 
	\atop
	\Leftrightarrow  B_0 \os B_1 \os ... \os B_i \os ... \os B_{i+n},...
	}
$$

\vspace{0.2cm}
\noindent results in \textbf{$F$-partial ordered set} $\langle\Phi,\leq\rangle$    with its Hasse digraph representation  looking like 
specific "cobweb"   image  (for cobweb posets portraits   see [1] and references therein  and see also cobwebs in action on\\
 $http://www.faces-of-nature.art.pl/cobwebposets.html$). 

\vspace{0.4cm}

\noindent\textbf{2.3  Partial order  $\leq$ .}

\vspace{0.1cm}
\noindent The partial order relation $\leq$ in the set of all points-vertices is determined  uniquely by the  above equivalent $F$- chains.  Let  $x,y \in  \Phi=\bigcup_{k=0}^n\Phi_k$ and let  $k, i  = 0,1,2,...$.   Then 

\begin{equation}\label{eq:leq}
	x\leq y \Leftrightarrow \forall_{x\in\Phi} : x\leq x \vee \Phi_i\ni x < y \in \Phi_{i+k}\ iff\ x(R_i\copyright...\copyright R_{i+k-1})y
\end{equation}

\noindent where  "$\copyright$"  stays for [Boolean] composition of binary relations.

\vspace{0.2cm}
\noindent \textbf{Relation  ($\leq$)is defined equivalently }: 

\vspace{0.2cm}
\noindent 
$ x  \leq  y$  in $(\Phi,\leq)$  iff either  $x=y$ or  there exist a directed path from $x$ to $y; x,y \in \Phi$.

\vspace{0.2cm}
\noindent Let now $R_k = \Phi_k\times\Phi_{k+1},  k \in N \cup\{0\}$. For "historical" reasons \cite{1}  we shall call such partial ordered set $\Pi = \langle\Phi,\leq\rangle$  the \textbf{cobweb poset} as  theirs  Hasse digraph representation  looks like specific "cobweb"   image.

\vspace{0.4cm}
\noindent \textbf{2.4. The natural join $\os$ operation and the natural join of  matrices satisfying the natural join condition}

\noindent We define here the adjacency matrices representation of the natural join $\os$ operation.

\vspace{0.1cm}
\noindent The adjacency matrix $\mathbf{A}$ of a bipartite graph with \textbf{biadjacency} = reduced adjacency  matrix $\mathbf{B}$ is given by

$$
	\mathbf{A} = \left(
	\begin{array}{cc}
	0 & \mathbf{B} \\
	\mathbf{B}^T & 0 \\
	\end{array}
	\right).
$$

\begin{defn}
The adjacency matrix $\mathbf{A}[D]$  of a bipartite \textbf{digraph} $D(R)= ( P\cup L ,  E \subseteq P\times L)$   with biadjacency matrix $\mathbf{B}$ is given by [1]

$$
	\mathbf{A}[D] = \left(
	\begin{array}{cc}
	0_{k,k} & \mathbf{B}(k\times m) \\
	0_{m,k} & 0_{m,m} \\
	\end{array}
	\right).
$$
where  $k = | P |$,   $m = | L |$.
\end{defn}

\vspace{0.2cm}

\noindent \textcolor{red}{\textbf{Note}}:  biadjacency and  cover relation $\prec\cdot$ matrix for bipartite digraphs coincide. By extension - we shall call  cover relation $\prec\cdot$ matrix also the biadjacency matrix.

\begin{conven}
 $S \copyright R$ = composition of  binary relations  $S$  and  $R  \Leftrightarrow   \mathbf{B}_{R\copyright S} =  \mathbf{B}_R \copyright \mathbf{B}_S$ where    ( $|V|= k , |W|= m$ )   $\mathbf{B}_R (k \times m) \equiv \mathbf{B}_R .$  
\end{conven}

\noindent $\mathbf{B}_R$ is the $(k \times m)$ \textbf{biadjacency} [or another name:  \textbf{reduced}   adjacency]  matrix  of the bipartite relations' $R$  digraph $B(R)$  and $\copyright$ apart from  relations composition  denotes also  Boolean multiplication of these rectangular biadjacency  Boolean matrices  $B_R , B_S$.    What is their form?   The answer is in the block structure of  the  standard square $(n \times n)$ adjacency matrix $A[D(R)];  n = k +m$ .  The form of  standard square adjacency matrix $A[G(R)]$ of bipartite digraph  $D(R)$ has the following  apparently  recognizable block reduced structure:  [ $O_{s\times s}$ stays for $(k \times m)$ zero matrix ]

$$
	\mathbf{A}[D(R)] = \left[
	\begin{array}{ll}
		O_{k\times k} & \mathbf{A}_R(k\times m) \\
		O_{m\times k} & O_{m\times m}
	\end{array}
	\right]
$$

\noindent Let $D(S) = (W(S)\cup T(S),E(S))$; $W\cap T = \emptyset$, $E (S)  \subseteq W\times T;$ ($|W|= m, |T|= s$); hence

$$
	\mathbf{A}[D(S)] = \left[
	\begin{array}{ll}
		O_{m\times m} & \mathbf{A}_S(m\times s) \\
		O_{s\times m} & O_{s\times s}
	\end{array}
	\right]
$$

\begin{defn} (natural join condition)
The ordered pair of matrices   $\langle \mathbf{A_1}, \mathbf{A_2} \rangle$ is said to satisfy the natural join condition iff  they  have the block structure  of     $\mathbf{A}[D(R)]$  and  $\mathbf{A}[D(S)]$  as above  i.e. iff  they might be identified accordingly : $\mathbf{A_1} =  \mathbf{A}[D(R)]$   and  $\mathbf{A_2} =  \mathbf{A}[D(S)]$.   
\end{defn}

\noindent Correspondingly if  two given digraphs $G_1$ and $G_2$  are such that their adjacency matrices  $\mathbf{A_1} = \mathbf{A}[G_1]$  and  $\mathbf{A_2} =\mathbf{A}[G_2]$ do satisfy the natural join condition we shall say  that $G_1$ and $G_2$  satisfy the natural join condition.
 For matrices satisfying the natural join condition one may define what follows.

\vspace{0.4cm}
\noindent  First we  define the \textbf{Boolean reduced}  or \textbf{natural join  composition} $\cs$   and secondly the natural join $\os$ of adjacent matrices  satisfying the natural join condition. 

\begin{defn} ($\cs$ composition)

$$
	\mathbf{A}[D(R\copyright S)] =: \mathbf{A}[D(R)] \cs \mathbf{A}[D(S)] =  \left[
	\begin{array}{ll}
		O_{k\times k} & \mathbf{A}_{R\copyright S}(k\times s) \\
		O_{s\times k} & O_{s\times s}
	\end{array}
	\right]
$$

\noindent where $\mathbf{A}_{R\copyright S}(k\times s) = \mathbf{A}_R(k\times m) \copyright \mathbf{A}_S(m\times s)$.
\end{defn}

\noindent according  to the scheme: 
$$
	[(k+m) \times (k + m )]  \cs [(m + s) \times (m + s)]  =  [(k+ s) \times (k+ s)] .
$$

\vspace{0.4cm}
\noindent \textbf{Comment 5.}
\noindent The adequate projection makes out the intermediate, joint in  common  $\mathrm{dom}(S) = \mathrm{rang}(R)=W$ , $|W|= m$.

\vspace{0.4cm}
\noindent The above Boolean reduced composition $\cs$  of adjacent matrices technically reduces then to the calculation of just Boolean product of  the  \textbf{reduced}  rectangular  adjacency matrices  of the bipartite relations` graphs.

\vspace{0.2cm}
\noindent We  are however  now in need of the Boolean natural  join product  $\os$  of adjacent matrices  already announced at the beginning of this presentation. Let us now define it.

\vspace{0.4cm}
\noindent As for  the \textbf{natural join} notion we aim at the morphism correspondence:
$$
	S \os R  \Leftrightarrow    M_{S\os R}  =  M_R \os M_S
$$

\noindent where $S \os R$ = natural  join of  binary relations  $S$  and  $R$  while
$M_{S\os R} =  M_R \os M_S$ = natural  join of  standard square adjacency matrices  
(with   customary convention: $M[G(R)]  \equiv  M_R$  adapted). Attention:   recall here that  the natural join of the above binary relations  $R \os S$  is  the ternary relation - and thus one results in $k$-ary relations  if with more factors undergo the $\os$ product.  As a matter of fact \textbf{ $\os$ operates on multi-ary relations according to the scheme:}   

$$
	(n+k)_{ary} \os (k+m)_{ary}  =   (n+ k +m)_{ary} .
$$

\noindent For example: $(1+1)_{ary} \os (1+1)_{ary} = (1+ 1 +1)_{ary}, binary \os binary = ternary$. 

\vspace{0.4cm}
\noindent Technically - the natural join of the $k$-ary  and  $n$-ary relations  is  defined accordingly the same way via  $\os$ natural join product of adjacency matrices - the adjacency matrices of   these relations' Hasse digraphs.

\vspace{0.2cm}
\noindent With the notation established above we finally define the natural join  $\os$  of two adjacency matrices as  follows:

\begin{defn} [natural join $\os$ of  matrices]. 

$$
 A[D(R \os S)]  =:  A[D(R)]  \os A[D(S)] =                                                
$$

$$
	= \left[
	\begin{array}{ll}
		O_{k\times k} & A_R(k\times m) \\
		O_{m\times k} & O_{m\times m}
	\end{array}
	\right] 
	\os
	\left[
	\begin{array}{ll}
		O_{m\times m} & A_S(m\times s) \\
		O_{s\times m} & O_{s\times s}
	\end{array}
	\right] =
$$
$$
	=\left[
	\begin{array}{lll}
		O_{k\times k} & A_R(k\times m) & O_{k\times s}\\
		O_{m\times k} & O_{m\times m}  & A_S(m\times s) \\
		O_{s\times k} & O_{s\times m}  & O_{s\times s}
	\end{array} 
	\right]
$$
\end{defn}

\vspace{0.2cm}
\noindent \textbf{Comment 6}.  The adequate projection used in natural join operation lefts one copy of the joint in  common "intermediate" submatrix $O_{m\times m}$ and consequently lefts one copy  of  "intermediate" joint in  common  $m$ according  to the scheme:
$$
	[(k+m) \times (k + m )]  \os  [(m + s) \times (m + s)]  =  [(k+ m + s) \times (k+ m + s)]  .
$$

\vspace{0.4cm}
\noindent \textbf{2.5.  The \textit{biadjacency i.e cover relation} $\prec\cdot$   matrices of  the  natural join of adjacency matrices}.

\noindent Denote  with  $B(A)$ the   biadjacency i.e cover relation $\prec\cdot$  matrix of the adjacency matrix $A$.\\
\textcolor{red}{\textbf{Note}}:  biadjacency and  cover relation $\prec\cdot$ matrix for bipartite digraphs coincide. By extension - we shall call  cover relation $\prec\cdot$ matrix the biadjacency matrix too, as for any graded digraph with more than one level we might represent it as a partition into
two independent sets, though it is more natural to see it as the natural join of the sequence of bipartite digraphs

\vspace{0.2cm}
\noindent Let  $A(G)$ denotes the adjacency matrix of the digraph $G$  , for example a di-biclique relation digraph.   Let  $A(G_k)$, $k= 0,1,2,...$ be the sequence adjacency matrices of the sequence $G_k, k=0,1,2,...$ of digraphs.  Let us identify  $B(A)\equiv B(G)$ as a convention.

\begin{defn} [digraphs natural join]
Let  digraphs $G_1$ and $G_2$   satisfy the natural join condition.  Let us make then  the identification  $A(G_1 \os G_2)  \equiv A_1 \os A_2$  as definition.  The digraph  $G_1 \os G_2$  is called the digraphs natural join of  digraphs  $G_1$ and $G_2$. Note that the order is essential. 
\end{defn}

\vspace{0.2cm}
\noindent We observe at once what follows. 

\begin{observen}
$$
	B (G_1 \os G_2) \equiv B (A_1 \os A_2) =  B(A_1)\oplus B(A_2) \equiv B (G_1)\oplus B(G_2)
$$
\end{observen}

\vspace{0.2cm}
\noindent \textbf{Comment 7.} The Observation 1 justifies the notation  $\os$ for the natural join
of relations digraphs and  equivalently for  the natural join of their adjacency matrices and 
equivalently for the natural join of relations that these are faithful representatives of.
\noindent Recall: $B(A)$ is the   biadjacency i.e cover relation $\prec\cdot$  matrix of the adjacency matrix $A$.\\
\textcolor{red}{\textbf{Note}}:  biadjacency and  cover relation $\prec\cdot$  matrix for bipartite digraphs coincide. Recall that by extension - we call  cover relation $\prec\cdot$ matrix the biadjacency matrix too.

\vspace{0.2cm}
\noindent As a consequence we have.

\begin{observen}
$$
	B\left(\os_{i=1}^n G_i \right) \equiv  B [\os_{i=1}^n  A(G_i)]  =  \oplus_{i=1}^n  B[A(G_i) ]  \equiv  \mathrm{diag} (B_1 , B_2 , ..., B_n) = 
$$
$$
	= \left[ \begin{array}{lllll}
	B_1 \\
	& B_2 \\
	& & B_3 \\
	& ... & ... & ...\\
	& & & & B_n
	\end{array} \right],
$$

\noindent or equivalently 

$$
	\kappa \equiv \chi \left(\prec\cdot \right) = \chi \left(\os_{i=1}^n \prec\cdot_i \right)  \equiv 
$$

$$
	\equiv \left[ \begin{array}{llllll}
	0 & B_1 \\
	& 0 & B_2 \\
	& & 0 & B_3 \\
	& ... & ... & ...\\
	& & & & 0 & B_n\\
	 & & & & & 0 

	\end{array} \right],
$$

\noindent $n \in N \cup \{\infty\}$ 
\end{observen}

\vspace{0.4cm}
\noindent \textbf{2.6. The formula  of zeta matrix for graded posets  with the finite set of minimal elements i.e for $F$-graded posets}

\vspace{0.2cm}
\noindent Any graded poset  with the finite set of minimal elements is an $F$- sequence denominated sub-poset of  its corresponding cobweb poset.
\noindent The Observation 2 supplies the simple recipe for the biadjacency (reduced adjacency) matrix of  Hasse digraph coding  any given  graded poset  with the finite set of minimal elements. The recipe for zeta matrix is then standard. We illustrate this by the source example; the source example as the adjacency  matrices  i.e  zeta matrices of any given  graded poset  with the finite set of minimal elements are sub-matrices of their corresponding cobweb posets and as such have the same block matrix structure. 
 
\vspace{0.2cm}
\noindent The explicit  expression for zeta matrix $\zeta_F$ of \textbf{\textcolor{blue}{cobweb posets}}  via known blocks of zeros and ones for arbitrary natural numbers valued $F$- sequence  was given in [1]  due to more than  mnemonic  efficiency  of the up-side-down notation being applied (see [1] and references therein). With this notation inspired by Gauss  and replacing  $k$ - natural numbers with   "$k_F$"  numbers one gets 
$$
	\mathbf{A}_F = \left[\begin{array}{llllll}
	0_{1_F\times 1_F} & I(1_F \times 2_F) & 0_{1_F \times \infty} \\
	0_{2_F\times 1_F} & 0_{2_F\times 2_F} & I(2_F \times 3_F) & 0_{2_F \times \infty} \\
	0_{3_F\times 1_F} & 0_{3_F\times 2_F} & 0_{3_F\times 3_F} & I(3_F \times 4_F) & 0_{3_F \times \infty} \\
	0_{4_F\times 1_F} & 0_{4_F\times 2_F} & 0_{4_F\times 3_F} & 0_{4_F\times 4_F} & I(4_F \times 5_F) & 0_{4_F \times \infty} \\
	... & etc & ... & and\ so\ on & ...
	\end{array}\right]
$$

\noindent and

$$
	\zeta_F = exp_\copyright[\mathbf{A}_F] \equiv (1 - \mathbf{A}_F)^{-1\copyright} \equiv I_{\infty\times\infty} + \mathbf{A}_F + \mathbf{A}_F^{\copyright 2} + ... =
$$
$$
	= \left[\begin{array}{lllll}
	I_{1_F\times 1_F} & I(1_F\times\infty) \\
	O_{2_F\times 1_F} & I_{2_F\times 2_F} & I(2_F\times\infty) \\
	O_{3_F\times 1_F} & O_{3_F\times 2_F} & I_{3_F\times 3_F} & I(3_F\times\infty) \\
	O_{4_F\times 1_F} & O_{4_F\times 2_F} & O_{4_F\times 3_F} & I_{4_F\times 4_F} & I(4_F\times\infty) \\
	... & etc & ... & and\ so\ on & ...
	\end{array}\right]
$$

\noindent where  $I (s\times k)$  stays for $(s\times k)$  matrix  of  ones  i.e.  $[ I (s\times k) ]_{ij} = 1$;  $1 \leq i \leq  s,  1\leq j  \leq k.$  and  $n \in N \cup \{\infty\}$
\vspace{0.2cm}

\begin{observen}
Let us denote by  $\langle\Phi_k\to\Phi_{k+1}\rangle$ (see the authors papers quoted) the di-bicliques  denominated by subsequent levels $\Phi_k, \Phi_{k+1}$ of the graded  $F$-poset $P(D) = (\Phi, \leq)$  i.e. levels $\Phi_k , \Phi_{k+1}$ of  its cover relation graded digraph  $D = (\Phi,\prec\!\!\cdot$)  [Hasse diagram].   Then

$$
	B\left(\os_{k=1}^n \langle\Phi_k\to\Phi_{k+1}\rangle \right) = \mathrm{diag}(I_1,I_2,...,I_n) = 
$$
$$
	= \left[ \begin{array}{lllll}
	I(1_F\times 2_F) \\
	& I(2_F\times 3_F) \\
	& & I(3_F\times 4_F) \\
	& & ... \\
	& & & & I(n_F\times (n+1)_F)
	\end{array} \right]
$$

\noindent where $I_k \equiv I(k_F \times (k+1)_F)$, $k = 1,...,n$ and where  - recall - $I (s\times k)$  stays for $(s\times k)$  matrix  of  ones  i.e.  $[ I (s\times k) ]_{ij} = 1$;  $1 \leq i \leq  s,  1\leq j  \leq k.$  and  $n \in N \cup \{\infty\}$.  
\end{observen}

\vspace{0.2cm}

\noindent The recipe for any $F$-denominated i.e.  the recipe for \textbf{\textcolor{red}{any  graded poset}} with a finite minimal elements set is supplied via the following observation.

\vspace{0.2cm}

\begin{observen}
Consider bigraphs'  chain \textbf{obtained from }the above di-biqliqes' chain via  deleting or no  arcs making thus [if deleting arcs] some or all of the di-bicliques $ \langle\Phi_k\to\Phi_{k+1}\rangle$  not di-biqliques; denote  them as  $G_k$. Let $B_k = B(G_k)$ denotes their biadjacency matrices correspondingly.  Then for any such $F$-denominated chain [hence any chain ] of bipartite digraphs  $G_k$  the general formula is:

$$
 B\left( \os_{i=1}^n G_i \right) \equiv  B [\os_{i=1}^n  A(G_i)] =  \oplus_{i=1}^n  B[A(G_i) ]  \equiv  \mathrm{diag} (B_1 , B_2 , ..., B_n) =
$$
$$
	= \left[ \begin{array}{lllll}
	B_1 \\
	& B_2 \\
	& & B_3 \\
	& & ... \\
	& & & & B_n
	\end{array} \right]
$$

\noindent $n \in N \cup \{\infty\}$.
\end{observen}

\vspace{0.2cm}

\noindent\textbf{Comment 8}  Note the notation identification: $	\zeta_F = exp_\copyright[\mathbf{A}_F] \equiv (1 - \mathbf{A}_F)^{-1\copyright}.$
\noindent Note that  $n!=1 mod 2$. Colligate also aside (?) reason: 
 $$lim_{q\rightarrow 1} exp_q  = exp$$ 
 
\noindent while
   $$lim_{q\rightarrow 0}exp_q[x] = (1-x)^{-1}. $$
    
\noindent Consult the Remark in [1] on the cases: Boolean  poset $2^N$ and the "Ferrand-Zeckendorf"  poset of finite subsets of $N$  without two consecutive elements. 

\vspace{0.2cm}

\begin{observen}
 The $F$-poset $P(G) = (\Phi, \leq)$  or equivalent to say:   its cover relation graded digraph $G = (\Phi,\prec\!\!\cdot) = \os_{k=0}^m G_k$ textcolor{blue}{\textbf{is of Ferrers dimension one}} \textbf{iff}  in the process of deleting arcs from  the cobweb poset Hasse diagram  $D = (\Phi,\prec\!\!\cdot)$ = $\os_{k=0}^n  \langle\Phi_k\to\Phi_{k+1}\rangle $   does not produces  $2\times 2$  permutation submatrices in any of  bigraphs $G_k$  biadjacency  matrices $B_k= B (G_k)$,  $k = 0,...,n$,  $n \in N \cup \left\{\infty \right\}$.
\end{observen}

\vspace{0.2cm}

\section{Cobweb posets and differential posets of Stanley  [3,4]}

\vspace{0.1cm}

\noindent \textbf{3.1. Preliminaries}

\vspace{0.1cm}

\noindent The graded  digraph $G$ (graph) is utterly denominated by its sequence of bipartite digraphs (graphs) that every two consecutive levels of $G$ do constitute. 

\noindent The complete  graded  digraph  $D$ is utterly denominated by its sequence of  complete bipartite digraphs - di-bicliques [1] that every two consecutive  levels of $D$  do constitute  ( see KoDAGs in [1] and references therein, consult the Example 2.7.2 in [5]). 

\noindent Because of their appearance, an at a first glance outlook - these complete  graded  digraph $D$ associated posets  where called   cobweb posets [ digraphs $D$ are identified with Hasse diagram of cobweb posets  see KoDAGs in [1,2] and references therein).

\vspace{0.1cm}

\noindent \textbf{Comment 9.}

\noindent The  appearance of "almost complete"  graded digraph (subgraphs of KoDAGs)  is  tremendously prevailing. These look like  - the hoary tree with silver cobweb threads.  The other extreme to the complete in such a  picture of a tree with cobweb (KoDAG) is  an also beautiful melancholic  bare rooted directed  tree graph - void of this spider's web  hoary tunicate and droplets (loops). 
\vspace{0.1cm}

\noindent Here are come some examples listed with the convention that graphs  become digraphs with all arcs directed upward (in the direction of increasing rank)  or dually - downwards  (in the direction of decreasing rank). If so has been done these become graded DAGs.  
\noindent First let us establish-recall for clarity that the directed tree (all arcs directed away from its root)  is a digraph which becomes a tree if directions on the edges are ignored. 
Colligate with an arborescence.  Naturally every arborescence is a directed acyclic graph.
\vspace{0.1cm}

\noindent Now come examples of sub-cobweb posets digraphs (hence DAGs) [5,3,4].

\noindent 0.  The binary tree digraph   (directed Tree)\\
\noindent 1.  The Fibonacci digraph   ([1],  directed Tree)\\
\noindent 2.  The Young graph \\                                                                  
\noindent 3.  The Young-Fibonacci graph\\                                                   
\noindent 4.  The Young-Fibonacci insertion graph\\                                     
\noindent 5.  The 2-dimensional Pascal graph\\                                           
\noindent 6.  The lattice of binary trees graph\\                                           
\noindent 7.  The lattice of Bracket  tree  graph\\                                         
\noindent 8.  The  $Fan_k$  graph  [5] (Amalgamate $k$ disjoint infinite chains by gluing their roots ("`zeros"')\\
\noindent 9.  The special complete graded graph [5] ( the set of vertices the same as in 8., connect any two vertices from consecutive 
levels  $\Phi_k$  and $\Phi_{k+1}$  with an upward directed arc) \\
\noindent 10. The $F$ - denominated (hence any) complete graded digraph with finite minimal elements set (KoDAG in [1,2] and references therein).

\vspace{0.3cm}

\noindent \textbf{3.2. Cobweb posets and differential posets}

\vspace{0.1cm}

\noindent The class of posets known as differential posets were  first introduced and studied by
Stanley in [3] then [4].  A differential poset from [3] is a partially ordered set with raising and lowering operators
$U$ and $D$ which satisfy the commutation relation $DU-UD = rI $ for some integer $r>0$. 
Generalizations of this class of posets were studied by Stanley [4] and Fomin [5]. A number of examples of generalized differential posets are given in these [3,4,5] papers. Another example, a poset of rooted unlabelled trees, was  introduced by Hoffman [6]. 

\vspace{0.2cm}
\noindent Let us consider at first the case $r=1$  which we shall call GHW case for the reasons to become apparent soon.

\vspace{0.1cm}
\noindent In this GHW  $r=1$ case the operators $U$ and $D$  are defined correspondingly, ($x,y,z \in \Phi$):

\begin{defn}

$$ Dx =  {\displaystyle \sum_{y \prec\cdot x}}y  $$

$$ Ux =  {\displaystyle \sum_{x\prec\cdot y}}y  $$
extended  by linearity  to (say complex, or...) linear space  $C[\Phi]$.

\end {defn}

\vspace{0.1cm}
\noindent From Theorem 2.2. in [3] we then have ( consult also [7]) that GHW commutation relation $DU-UD = I $ holds iff 
$P(D) = \left\langle \Phi, \leq \right\rangle  $ is differential ($r=1$) poset. Out of this one infers inductively [8.9,10,11,7] 
what follows.

\begin{observen}

$$DU^n = nU^{n-1} + U^nD,$$
for $n\in N.$

\end{observen}

\vspace{0.1cm}

\begin{observen}
\noindent Cobweb posets Hasse digraphs from 10. above are examples of  $\stackrel{\rightarrow}{q}$, $\stackrel{\rightarrow}{r}$ - differential posets (see [5] for all $q_n = 1$)  and might serve for  more general structures  (due to Fomin [5] )  called  Dual graded Graphs.
\end{observen}

\noindent \textbf{Check.} Indeed. In [5]  $U_n$  and $D_n$  are defined as restrictions of  $U$  and  $D$  onto $C[\Phi_k]$"`homogeneous"' subspaces of $C[\Phi]$, 
$k=0,1,2,...$. Then 

$$ D_{n+1} = q_n U_{n-1} D_n + r_n I_n, \:  n\in N,$$

\vspace{0.1cm}

\noindent where  $r_0= 1_F, q_0 = 0$  for $n = 0$ and $q_n = \frac{(n+1)_F}{(n-1)_F}, r_n = 0$ for $n>0$. 

\noindent For  Fomin examples 8. and 9.  above  (see 2.7.2. in [5])  the $(1.4.11)$ i.e. the $(1.4.10)$ with all $q_n = 1$  from  [5] holds  
for   $r_0 = t\equiv 0_F,  r_1  = r_2 = r_3...  =  0$

\vspace{0.1cm}

\noindent The KoDAGs graded digraphs example 10. in view of Observation 7  and Observation 6 becomes the motivating example of the following description
appealing to the corresponding R. Stanley definition from [3]. One expect this description  to be  efficient and tangible while the combinatorics of path counting is concerned also via  colligation  with  quantum models - for example.

\begin{defn}

\noindent Let $F$ be such that $0_F=1$. The locally finite graded poset $\Pi = \left\langle \Phi, \leq \right\rangle$ is then said to be $F$-differential poset iff\\

\noindent \textbf{1.}  if  $x,y \in \Phi, x \neq y$ and there are $k_F$ elements of $\Phi$ which cover both of them then there are exactly $k_F$ elements of $\Phi$ which are covered by both of them. 

\noindent \textbf{2.} If  $x \in \Phi$ covers exactly $k_F$ elements from $\Phi$  then this very $x$ is covered exactly by $(k+1)_F$ elements from $\Phi$ 
\end{defn}

\noindent As in  $F=N$ case  (i.e  for all  $k = 1,2,...$ we have  $k_F \equiv k_N \equiv k $)  the above conditions determine the number of elements and cover relations among them in the up to fourth rank [7]. Above the fourth rank ($k>4$) cover relations and number of elements in levels $\Phi_k$ vary.
\noindent And what does happen when the requirement $0_F=1$ is relaxed?

\vspace{0.2cm} 

\noindent Let us introduce the direct sum of certain projection on $\Phi_n$ operators and denote this operator  with the symbol $\delta_F$.
Let $x_n$  denotes homogeneous element of $\Phi$   i.e  $x \in \phi_n$.  Then  $\delta_F$  is specified as follows.

\begin{defn} 

$\delta_F$ = $diag(1_F - 0_F, 2_F - 1_F, 3_F -2_F,..., (n+1_F - n_F, ...)\equiv diag(\delta_0, \delta_1,\delta_2,...,\delta_n, ...)$ i.e.   $ \delta_F (x_n) = \delta_n x_n . $

\end{defn}

\vspace{0.2cm} 

\noindent A straightforward verifying (see [7]) leads us to thus confirmed conclusion below; (note that $\delta_N = I$).

\begin{observen}

\noindent $$DU - UD = \delta_F  $$

\end{observen}

\vspace{0.2cm} 

\noindent Out of this (as in [10,11])  one infers inductively  
what follows.

\begin{observen}

$$DU^n = n\delta_F U^{n-1} + U^nD,$$
for $n\in N.$
\end{observen}

\vspace{0.2cm} 

\noindent Since $Dx_0= 0$ we have  $DU = n \delta_F U^{n-1}x_0$ which means that $D$ is 
representative of the Markowsky general linear operator i.e. a derivative from extended umbral calculus (see [10,11,12] and plenty of references therein). 

\vspace{0.2cm} 

\noindent The $F$-denominated cobweb posets are in a sense a canonical example of  $F$- differential posets. More on that is expected soon.

\vspace{0.2cm} 

\noindent \textbf{Comment 10}  (  miscellaneous (aside?) final  remark)\\
\noindent The  ingenious ideas of differential and dual graded posets that we owe to  Stanley [3,4]  and  Fomin [5] bring together combinatorics, representation theory, topology, geometry  and many more specific branches of mathematics and mathematical physics thanks to intrinsic ingredient of these mathematical descriptions which is the Graves - Heisenberg - Weyl  (GHW)  algebra usually attributed to Heisenberg by physicists and  to Herman Weyl by mathematicians  and sometimes  to both of them (see: [3] for Weyl, [5] for Heisenberd and then [8] and [9] ; for  GHW  see  [10-12]  then  note the content and context of  [13,14] ).
\noindent As noticed by the author of  [9]  the formula  

$$ [f(a),b] = cf'(a)$$
                                                        
\noindent where
$$[a,b]=c,  [a,c]=[b,c] = 0 $$
                                                      
\noindent pertains to Charles Graves from Dublin [8]. Then it was re-discovered by Paul Dirac and others in the next century.\\
Let us then note that the picture that emerges in [11-12] discloses the fact that any umbral representation of finite (extended) operator calculus or equivalently - any umbral representation of GHW algebra makes up an example of the algebraization of the analysis with generalized differential operators  acting on the algebra of polynomials or other algebras as for example formal series algebras.

\vspace{0.4cm}

\noindent \textbf{Bibliography remark.} On  Umbra Difference Calculus references streams see [15] including references ad Comment 8  and all of that. On the history of cobweb poset $\zeta$ function formulas and also for the $\zeta_{-1}$ formula see [16].

\end{document}